# UTILITY MAXIMIZING ENTROPY AND THE SECOND LAW OF THERMODYNAMICS[1]

By Wojciech Słomczyński and Tomasz Zastawniak

*Jagiellonian University and University of Hull*

*Dedicated to Professor Andrzej Lasota on his 70th birthday*

Expected utility maximization problems in mathematical finance lead to a generalization of the classical definition of entropy. It is demonstrated that a necessary and sufficient condition for the second law of thermodynamics to operate is that any one of the generalized entropies should tend to its minimum value of zero.

**1. Introduction.** The maximization of expected logarithmic utility is well known to be related to the classical notion of Boltzmann–Gibbs entropy $H(f) = \int_\Omega f \ln f \, d\mu$, namely

$$H(f) = \sup_w \int_\Omega f \ln w \, d\mu$$

for any density $f$, the supremum being taken over all densities $w$ under the probability measure $\mu$ on $\Omega$. This is a consequence of the integrated Gibbs inequality $\int_\Omega f \ln f \, d\mu \geq \int_\Omega f \ln w \, d\mu$, valid for any densities $f$ and $w$ (see, e.g., [20, 22]).

Several authors, including Bismut [3], Pikovsky and Karatzas [23], Amendinger, Imkeller and Schweizer [1], Frittelli [9, 10], Bellini and Frittelli [2], Schachermayer [27], Kramkov and Schachermayer [17, 18], Cvitanić, Schachermayer and Wang [7], Delbaen et al. [8], Rouge and El Karoui [26], Goll and Rüschendorf [12] and others, developed duality methods in the context of semimartingale theory, and in recent years have applied them in mathematical finance to investigate

---

Received July 2002; revised June 2003.

[1]Supported in part by the British Council and the Polish State Committee for Scientific Research Grant WAR/341/204.

*AMS 2000 subject classifications.* Primary 37A50; secondary 94A17, 60F25, 91B16, 49N15.

*Key words and phrases.* Entropy, Markov operators, utility maximization, exactness, *H*-theorem, second law of thermodynamics.

$\sup_w \int_\Omega u(w) \, d\nu$ over all possible finite values $w \geq 0$ of self-financing trading strategies with fixed initial wealth $a$ for a general class of utility functions $u$, where $\nu$ is a probability measure that captures the true probabilities of possible market scenarios. According to [9, 10] and [2], in a wide class of arbitrage-free markets, there is a pricing measure $\mu$, called a minimax martingale measure, such that this supremum is equal to $\sup_w \int_\Omega w \, u(w) \, d\nu$ over all random variables $w \geq 0$ with $\int_\Omega w \, d\mu = a$; see also [[11]–[15]].

We take unit initial wealth $a = 1$. When $\mu \gg \nu$ with density $f = \frac{d\nu}{d\mu}$, the last supremum can be written as $\sup_w \int_\Omega f u(w) \, d\mu$ and is taken over all densities $w$ with respect to measure $\mu$. In particular, for the logarithmic utility $u(x) = \ln x$, $x > 0$, the supremum is equal to $H(f) = u(e^{H(f)})$, where

$$H(f) = \int_\Omega f \ln f \, d\mu$$

is the Boltzmann–Gibbs entropy of $f$. Moreover, in the case of isoelastic utility $u(x) = \frac{1}{\gamma} x^\gamma$, $x > 0$, where $\gamma \in (-\infty, 0) \cup (0, 1)$, it is not hard to verify that the supremum is equal to $u(e^{H_\alpha(f)})$, where

$$H_\alpha(f) = \frac{1}{\alpha - 1} \ln \int_\Omega f^\alpha \, d\mu$$

is the Rényi entropy of order $\alpha = (1 - \gamma)^{-1} \in (0, 1) \cup (1, \infty)$; see [25].

These observations suggest that for a large class of utility functions $u$, the functional $H_u(f)$, defined by

$$u(e^{H_u(f)}) = \sup_w \int_\Omega f u(w) \, d\mu$$

for any density $f$ under $\mu$, where the supremum is taken over all densities $w$ with respect to $\mu$, may share some general properties of the Boltzmann–Gibbs entropy $H(f)$ or the Rényi entropy $H_\alpha(f)$. We propose to call $H_u(f)$ the utility maximizing entropy or $u$-entropy; see Definition 4.

The class of utility functions considered here consists of all strictly concave, strictly increasing, continuously differentiable functions $u : (0, \infty) \to \mathbb{R}$ such that $\lim_{x \searrow 0} u'(x) = \infty$ and $\lim_{x \nearrow \infty} u'(x) = 0$, satisfying the asymptotic elasticity condition

$$AE(u) = \limsup_{x \nearrow \infty} \frac{x u'(x)}{u(x)} < 1$$

of Kramkov and Schachermayer [17]. The asymptotic elasticity condition is imposed to ensure that the supremum is realized for some density $w$; see Theorem 20.

In this paper we demonstrate that $H_u(f)$ plays a similar role in the thermodynamic equilibrium limit as the classical Boltzmann–Gibbs entropy $H(f)$.

The states of a thermodynamic system are identified with the densities $f$ on a phase space $\Omega$ equipped with measure $\mu$. The evolution of a state $f$ can be described in terms of the iterations $f, Pf, P^2f, \ldots$ of a Markov operator $P$, that is, a linear operator on $L^1(\mu)$ that transforms densities into densities. The existence of a stationary density $f = Pf$ corresponds to a state of thermodynamic equilibrium.

The *second law of thermodynamics* (in its strong form) stipulates the existence of only one state $f$ of thermodynamic equilibrium that is approached regardless of the initial state of the system and is associated with the minimum value zero of Boltzmann–Gibbs entropy $H(f)$. On a space of measure 1 this state must necessarily be given by the uniform density $f \equiv 1$. A necessary and sufficient condition for the second law to operate is that the Markov operator $P$ should be *exact*, that is, $P^nf \to 1$ in $L^1(\mu)$ as $n \to \infty$ for any density $f$; equivalently, the Boltzmann–Gibbs entropy should tend to its minimum value of zero, $H(P^nf) \searrow 0$ as $n \to \infty$ for any density $f$ such that $H(f) < \infty$; see [20] or [22].

The main result of this paper, Theorem 29, is that the Boltzmann–Gibbs entropy can be replaced by the $u$-entropy for any given utility function $u$ that satisfies the asymptotic elasticity condition. That is to say, $P^nf \to 1$ in $L^1(\mu)$ as $n \to \infty$ for any density $f$ if and only if $H_u(P^nf) \searrow 0$ as $n \to \infty$ for any density $f$ such that $H_u(f) < \infty$. In other words, $u$-entropy can play exactly the same role in the second law as the Boltzmann–Gibbs entropy. The results also extend to Markov semigroups; see Theorem 31.

The behavior of Boltzmann–Gibbs entropy under the action of a Markov operator has been studied by many authors.

The fact that the sequence $H(P^nf)$ is decreasing ($H$-theorem) can easily be derived from the Jensen inequality for Markov operators. The idea goes back at least as far as the early papers of Csiszár (see also [[20, 21, 28]]).

The implication

$$(1) \qquad H(f_n) \to 0 \quad \Longrightarrow \quad f_n \xrightarrow{L^1} 1,$$

which is true for an arbitrary sequence of densities $\{f_n\}_{n \in \mathbb{N}}$, follows immediately from the Pinsker–Kullback-Csiszár inequality: $\frac{1}{2}\|f - 1\|_{L^1}^2 \leq H(f)$ (see [[5, 19, 24]]; for another proof, see [21]). In fact, Łoskot and Rudnicki [21] proved this implication for a larger class of entropy-like quantities, Csiszár's $\eta$-entropies [5] $H_\eta(f) = \int_\Omega \eta(f)\,d\mu$, where $\eta \colon [0, \infty) \to \mathbb{R}$ is an arbitrary convex function such that $\eta(0) = 0$. The result was applied in [4] and [21] to analyze the stability of solutions of parabolic equations. The notion of $\eta$-entropy also covers the case of Rényi entropy of order $\alpha \in (0, 1) \cup (1, \infty)$. The proof of the implication in [20] applies to a sequence of the special form $f_n = P^nf$ ($n \in \mathbb{N}$) with $P$ a Markov operator and uses the Komornik–Lasota spectral decomposition theorem for Markov operators [20].

The reverse implication to (1) is not true in general; see [6] for a counterexample. For sequences of the form $f_n = P^n f$ $(n \in \mathbb{N})$ with $H(f) < \infty$ and $P$ a Markov operator, the implication

$$(2) \qquad\qquad f_n \xrightarrow{L^1} 1 \quad \Longrightarrow \quad H(f_n) \to 0$$

was proved in [20].

Our results in Theorems 27 and 26 generalize both implications (1) and (2) to the case of $u$-entropy. The proof of Theorem 27 rests on the data-reduction inequality technique invented by Csiszár. To prove Theorem 26 we derive a formula for $u$-entropy by convex duality methods, inspired by the work of Kramkov and Schachermayer [17, 18].

To conclude the introductory part, let us remark that mathematical finance has enjoyed unprecedented growth in recent years, not least because of considerable input from other disciplines, including physical sciences, in general, and thermodynamics, in particular. Here we see it returning the favor: The $u$-entropy introduced above, hinted upon in many recent works on expected utility maximization, turns out to play a major role in the second law of thermodynamics, on an equal footing with the classical notion of Boltzmann–Gibbs entropy.

### 1.1. *Notation.*
The following notational conventions are used throughout the paper:

1. Take $\infty \cdot 0 = 0$ and $-\infty \cdot 0 = 0$.
2. Take $\mathbb{R}^+ = (0, \infty)$.
3. Take $(\Omega, \Sigma, \mu)$ to be a probability space.
4. Take $D(\mu)$ to denote the set of all densities on $(\Omega, \Sigma, \mu)$,

$$D(\mu) = \left\{ w \in L^1(\mu) : w \geq 0 \text{ and } \int_\Omega w \, d\mu = 1 \right\}.$$

5. Take $f\mu$ to be the probability measure absolutely continuous with respect to $\mu$ with density $f \in D(\mu)$, that is, for any $A \in \Sigma$,

$$(f\mu)(A) = \int_A f \, d\mu.$$

6. Take $\|\cdot\|_\alpha$ to denote the norm in $L^\alpha(\mu)$ for any $\alpha \in [1, \infty]$, and a pseudonorm for any $\alpha \in (0, 1)$.

### 2. Entropy.
In this section we define $u$-entropy and establish its principal properties. In particular, in Theorem 20 we obtain a formula for $u$-entropy.

2.1. *Utility functions.* We begin by recalling the definitions and properties concerned with utility functions and convex analysis that are needed throughout this paper.

DEFINITION 1. Let $u : \mathbb{R}^+ \to \mathbb{R}$. We call $u$ a *utility function* whenever $u$ is strictly concave, strictly increasing, continuously differentiable and such that

$$u'(0) = \lim_{x \searrow 0} u'(x) = \infty, \qquad u'(\infty) = \lim_{x \nearrow \infty} u'(x) = 0.$$

We also use the notation

$$u(0) = \lim_{x \searrow 0} u(x), \qquad u(\infty) = \lim_{x \nearrow \infty} u(x).$$

PROPOSITION 1. *The function* $I = (u')^{-1} : \mathbb{R}^+ \to \mathbb{R}^+$ *is strictly decreasing and satisfies*

$$I(0) = \lim_{x \searrow 0} I(x) = \infty, \qquad I(\infty) = \lim_{x \nearrow \infty} I(x) = 0.$$

DEFINITION 2. Let $u : \mathbb{R}^+ \to \mathbb{R}$ be a utility function. The *convex dual* $u^* : \mathbb{R}^+ \to \mathbb{R}$ is defined by

$$(3) \qquad u^*(y) = \sup_{x > 0} (u(x) - yx)$$

for any $y \in \mathbb{R}^+$.

The following basic properties of convex functions and convex duals can be found in various books, for example, [16].

PROPOSITION 2. *Let* $u : \mathbb{R}^+ \to \mathbb{R}$ *be a utility function. Then:*

1. *The function* $u^*$ *is strictly convex, strictly decreasing and continuously differentiable.*
2. *The equalities* $u^*(0) = u(\infty)$, $u^*(\infty) = u(0)$, $(u^*)'(0) = -\infty$ *and* $(u^*)'(\infty) = 0$ *hold.*
3. *For any* $y \in \mathbb{R}^+$,

$$(4) \qquad u^*(y) = u(I(y)) - yI(y).$$

4. *For any* $x \in \mathbb{R}^+$,

$$(5) \qquad u(x) = \inf_{y > 0} (u^*(y) + xy).$$

5. *For any* $y \in \mathbb{R}^+$,

$$(u^*)'(y) = -I(y).$$

Lemma 3. *Let* $u : \mathbb{R}^+ \to \mathbb{R}$. *Then* $u$ *is concave if and only if there exist* $a_n > 0$, $b_n \in \mathbb{R}$ *for any* $n \in \mathbb{N}$ *such that* $u(x) = \inf\{a_n x + b_n : n \in \mathbb{N}\}$ *for every* $x > 0$.

Example 1. Let $\gamma \in (-\infty, 1)$. Define $u_\gamma : \mathbb{R}^+ \to \mathbb{R}$ by

$$u_\gamma(x) = \begin{cases} \dfrac{1}{\gamma} x^\gamma, & \text{for } x \in \mathbb{R}^+ \text{ and } \gamma \in (-\infty, 0) \cup (0, 1), \\ \ln x, & \text{for } x \in \mathbb{R}^+ \text{ and } \gamma = 0. \end{cases}$$

We call $u_\gamma$ the *isoelastic utility of order* $\gamma$ if $\gamma \neq 0$ and the *logarithmic utility* if $\gamma = 0$.

The following definition of asymptotic elasticity and its properties is due to Kramkov and Schachermayer [17].

Definition 3. Let $u : \mathbb{R}^+ \to \mathbb{R}$ be a utility function. Then we define the *asymptotic elasticity* of $u$ by

$$\mathrm{AE}(u) = \limsup_{x \nearrow \infty} \frac{x u'(x)}{u(x)}.$$

Note that $\mathrm{AE}(u_\gamma) = \gamma$ for $\gamma < 1$.

Proposition 4. *Let* $u : \mathbb{R}^+ \to \mathbb{R}$ *be a utility function. Then*

$$\mathrm{AE}(u) \in \begin{cases} [0, 1], & \text{if } u(\infty) = \infty, \\ \{0\}, & \text{if } 0 < u(\infty) < \infty, \\ [-\infty, 0], & \text{if } -\infty < u(\infty) \leq 0. \end{cases}$$

Proposition 5. *Let* $u : \mathbb{R}^+ \to \mathbb{R}$ *be a utility function, let* $a > 0$ *and let* $b \in \mathbb{R}$. *Then* $\widetilde{u} = au + b : \mathbb{R}^+ \to \mathbb{R}$ *is a utility function. If* $u(\infty), \widetilde{u}(\infty) > 0$, *then* $\mathrm{AE}(u) = \mathrm{AE}(\widetilde{u})$.

Corollary 6. *Let* $u : \mathbb{R}^+ \to \mathbb{R}$ *be a utility function and let* $\widetilde{u} = au + b$ *for* $a > 0$, $b \in \mathbb{R}$. *Then* $\mathrm{AE}(u) < 1$ *if and only if* $\mathrm{AE}(\widetilde{u}) < 1$.

Proposition 7. *Let* $u : \mathbb{R}^+ \to \mathbb{R}$ *be a utility function such that* $u(\infty) > 0$ *and* $\mathrm{AE}(u) < \gamma < 1$. *Then there is an* $x_0 > 0$ *such that* $0 < u(\lambda x) \leq \lambda^\gamma u(x)$ *for all* $\lambda \geq 1$ *and* $x \geq x_0$.

Proof. The following argument slightly simplifies the proof of Kramkov and Schachermayer [17]: From the definition of asymptotic elasticity, it follows that there exists $x_0 > 0$ such that $0 < \gamma u(x) - x u'(x)$ for any $x \geq x_0$. For such an $x$ we define a function $G_x : [1, \infty) \to \mathbb{R}$ by $G_x(\lambda) = \lambda^\gamma u(x) - u(\lambda x)$ for $\lambda \geq 1$. Then $G_x(1) = 0$ and $G'_x(1) = \gamma u(x) - xu'(x) > 0$. Moreover, $G'_x(\lambda) = \gamma \lambda^{\gamma-1} u(x) - xu'(\lambda x) = \frac{\gamma}{\lambda}(G_x(\lambda) + u(\lambda x) - \frac{\lambda}{\gamma} xu'(\lambda x)) > \frac{\gamma}{\lambda} G_x(\lambda)$ for $\lambda > 1$. Using the theory of differential inequalities, we can deduce that $G_x(\lambda) \geq 0$ for all $\lambda \geq 1$, which completes the proof. $\square$

2.2. *Definition and basic properties of $u$-entropy.* Throughout the rest of this paper $u : \mathbb{R}^+ \to \mathbb{R}$ denotes a utility function in the sense of Definition 1.

DEFINITION 4. For any $f \in D(\mu)$ we put

$$N_u(f) = \sup_{w \in \mathcal{A}(f)} \int_\Omega u(w) f \, d\mu,$$

where

$$\mathcal{A}(f) = \{w \in D(\mu) : u(w)^- \in L^1(f\mu)\}.$$

Note that $\int_\Omega u(w) f \, d\mu \in (-\infty, \infty]$ for each $w \in \mathcal{A}(f)$. Now we define

$$H_u(f) = \ln u^{-1}(N_u(f))$$

and call it the *utility maximizing entropy* or *$u$-entropy* of $f$.

PROPOSITION 8. *The following inequalities hold for any $f \in D(\mu)$:*

$$u(1) \leq N_u(f) \leq u(\infty),$$
$$0 \leq H_u(f) \leq \infty.$$

PROOF. Taking $w \equiv 1 \in \mathcal{A}(f)$, we obtain the lower bound. The upper bound follows immediately from the definition. $\square$

PROPOSITION 9. *For any $f \in D(\mu)$, the following conditions are equivalent:*

1. $N_u(f) < u(\infty)$;
2. $N_u(f) < \infty$;
3. $H_u(f) < \infty$.

*In particular, all three conditions are satisfied for any utility function $u$ such that $u(\infty) < \infty$.*

PROOF. The implications $1 \Rightarrow 3 \Rightarrow 2$ are obvious, as is $2 \Rightarrow 1$ when $u(\infty) = \infty$.

Let us prove $2 \Rightarrow 1$ when $u(\infty) < \infty$. In this case, take an $n \in \mathbb{N}$ such that $\mu(A) \leq \frac{1}{n} \Rightarrow (f\mu)(A) \leq \frac{1}{2}$ for each measurable set $A$. Such an $n$ exists because $f\mu$ is absolutely continuous with respect to $\mu$. Let $w \in \mathcal{A}(f)$. Then $\mu\{w \geq n\} \leq \frac{1}{n}$ because $\int_\Omega w \, d\mu = 1$. It follows that $(f\mu)\{w \geq n\} \leq \frac{1}{2}$. As a result,

$$\begin{aligned}
\int_\Omega u(w) f \, d\mu &= \int_{\{w < n\}} u(w) f \, d\mu + \int_{\{w \geq n\}} u(w) f \, d\mu \\
&\leq u(n)(1 - (f\mu)\{w \geq n\}) + u(\infty)(f\mu)\{w \geq n\} \\
&= u(n) + (u(\infty) - u(n))(f\mu)\{w \geq n\} \\
&\leq u(n) + \frac{u(\infty) - u(n)}{2} = \frac{u(\infty) + u(n)}{2}.
\end{aligned}$$

Hence

$$N_u(f) = \sup_{w \in \mathcal{A}(f)} \int_\Omega u(w) f \, d\mu \leq \frac{u(\infty) + u(n)}{2} < u(\infty),$$

as required. This also shows that all three conditions must be satisfied whenever $u(\infty) < \infty$, completing the proof. □

PROPOSITION 10. *For any $f \in D(\mu)$, the following conditions are equivalent:*

1. $H_u(f) = 0$;
2. $f = 1$ ($\mu$-a.e.).

PROOF. $2 \Rightarrow 1$. Let $f = 1$ $\mu$-a.e. Take any $w \in \mathcal{A}(f)$. By Jensen's inequality, $\int_\Omega u(w) \, d\mu \leq u(\int_\Omega w \, d\mu) = u(1)$. Hence $N_u(f) = u(1)$ and $H_u(f) = 0$ by Proposition 8.

$1 \Rightarrow 2$. Suppose that $\mu\{f = 1\} < 1$ and take $A = \{f > 1\}$. Since $f$ is a density under $\mu$, it follows that $0 < \mu(A) < 1$. Moreover, $(f\mu)(A) > \mu(A)$. For any $a \in (0, 1)$ we put

$$w_a = \frac{a}{\mu(A)} \mathbb{1}_A + \frac{1 - a}{\mu(A^\neg)} \mathbb{1}_{A^\neg} \in \mathcal{A}(f),$$

$$\varphi(a) = \int_\Omega u(w_a) f \, d\mu = u\left(\frac{a}{\mu(A)}\right)(f\mu)(A) + u\left(\frac{1 - a}{\mu(A^\neg)}\right)(f\mu)(A^\neg).$$

Clearly, $w_{\mu(A)} \equiv 1$ and $\varphi(\mu(A)) = u(1)$. Moreover,

$$\varphi'(\mu(A)) = u'(1)\frac{(f\mu)(A) - \mu(A)}{\mu(A)\mu(A^\neg)} > 0.$$

As a result, there is an $a \in (0, 1)$ such that $\varphi(a) > \varphi(\mu(A))$. Because $N_u(f) \geq \varphi(a)$ for any $a \in (0, 1)$, it follows that $N_u(f) > u(1)$ and $H_u(f) > 0$. □

PROPOSITION 11. *Let $a > 0$ and let $b \in \mathbb{R}$. If $u : \mathbb{R}^+ \to \mathbb{R}$ is a utility function, then $au + b$ is also a utility function, and for any $f \in D(\mu)$,*

$$N_{au+b}(f) = aN_u(f) + b,$$
$$H_{au+b}(f) = H_u(f).$$

These properties follow immediately from the definition.

REMARK 1. Using Proposition 11, in many arguments we can assume without loss of generality that $u(1) = 0$.

PROPOSITION 12. *Let $f \in D(\mu)$. We define*

$$\mathcal{A}_b(f) = \{w \in D(\mu) : u(w)^- \in L^1(f\mu), w \text{ is bounded}\},$$
$$\mathcal{A}_0(f) = \{w \in D(\mu) : u(w)^- \in L^1(f\mu), w\mu \ll f\mu\},$$
$$\mathcal{A}_{0b}(f) = \{w \in D(\mu) : u(w)^- \in L^1(f\mu), w \text{ is bounded}, w\mu \ll f\mu\}.$$

*Then*

$$N_u(f) = \sup_{w \in \mathcal{A}_b(f)} \int_\Omega u(w) f \, d\mu = \sup_{w \in \mathcal{A}_0(f)} \int_\Omega u(w) f \, d\mu$$
$$= \sup_{w \in \mathcal{A}_{0b}(f)} \int_\Omega u(w) f \, d\mu.$$

PROOF. Put $\mathcal{A}_i = \mathcal{A}_i(f)$, $N_i = \sup_{w \in \mathcal{A}_i} \int_\Omega u(w) f \, d\mu$ for $i = b, 0, 0b$, $N = N_u(f)$ and $\mathcal{A} = \mathcal{A}(f)$. Clearly, $\mathcal{A}_{0b} \subset \mathcal{A}_0 \subset \mathcal{A}$ and $\mathcal{A}_{0b} \subset \mathcal{A}_b \subset \mathcal{A}$. Hence $N_{0b} \leq N_0 \leq N$ and $N_{0b} \leq N_b \leq N$. We can assume without loss of generality that $u(1) = 0$.

We show that $N \leq N_b$. Let $w \in \mathcal{A}$. Define $w_n = w \mathbb{1}_{\{w < n\}} + a_n \mathbb{1}_{\{w \geq n\}}$, for $n = 1, 2, \dots$, where $a_n = \mu\{w \geq n\}^{-1} \int_{\{w \geq n\}} w \, d\mu \geq n$. Then $\int_\Omega w_n \, d\mu = 1$ and $w_n \to w$ as $n \to \infty$. Clearly, $0 \leq w_n \leq a_n$, so $w_n$ is bounded. Moreover, $\{w_n < 1\} = \{w < 1\}$. Hence $u(w_n)^- = u(w)^- \in L^1(f\mu)$. Thus $w_n \in \mathcal{A}_b$ for every $n \in \mathbb{N}$. Applying the Fatou lemma, we obtain $\int_\Omega u(w) f \, d\mu \leq \liminf_{n \to \infty} \int_\Omega u(w_n) f \, d\mu \leq N_b$. Consequently, $N \leq N_b$.

Next, we show that $N_b \leq N_{0b}$. Let $w \in \mathcal{A}_b$. If $\int_{\{f > 0\}} w \, d\mu > 0$, then we take

$$\widetilde{w} = \frac{\mathbb{1}_{\{f > 0\}} w}{\int_{\{f > 0\}} w \, d\mu}.$$

Clearly, $\widetilde{w} \in D(\mu)$ is bounded and $\widetilde{w}\mu \ll f\mu$. Since $w \leq \widetilde{w}$ on $\{f > 0\}$, it follows that $\int_\Omega u(w)^- f \, d\mu = \int_{\{f > 0\}} u(w)^- f \, d\mu \geq \int_{\{f > 0\}} u(\widetilde{w})^- f \, d\mu \geq 0$, so $u(\widetilde{w})^- \in L^1(f\mu)$. As a result, $\widetilde{w} \in \mathcal{A}_{0b}$. Observe that $\int_\Omega u(w)f \, d\mu = \int_{\{f>0\}} u(w)f \, d\mu \le \int_{\{f>0\}} u(\widetilde{w})f \, d\mu \le N_{0b}$. If, on the other hand, $\int_{\{f>0\}} w \, d\mu = 0$, then we take

$$\widehat{w} = \frac{\mathbb{1}_{\{f>0\}}}{\mu\{f>0\}},$$

which clearly also belongs to $\mathcal{A}_{0b}$. Moreover, $\int_\Omega u(w)f \, d\mu = u(0) < u(1) \le \int_\Omega u(\widehat{w})f \, d\mu \le N_{0b}$ because $w = 0$ and $\widehat{w} \ge 1$ on $\{f > 0\}$ ($\mu$-a.e.). As a consequence, $N_b \le N_{0b}$.  $\square$

PROPOSITION 13.  *Let $f_1, f_2 \in D(\mu)$ and let $a \in [0, 1]$. Then*

$$N_u(af_1 + (1-a)f_2) \le aN_u(f_1) + (1-a)N_u(f_2).$$

PROOF.  Put $f = af_1 + (1-a)f_2$. First observe that $\int_\Omega u^-(w)f \, d\mu = a \int_\Omega u^-(w)f_1 \, d\mu + (1-a) \int_\Omega u^-(w)f_2 \, d\mu$ for any $w \in D(\mu)$, so $\mathcal{A}(f) = \mathcal{A}(f_1) \cap \mathcal{A}(f_2)$. Hence

$$\begin{aligned}
N_u(f) &= \sup_{w \in \mathcal{A}(f)} \int_\Omega u(w)(af_1 + (1-a)f_2) \, d\mu \\
&\le a \sup_{w \in \mathcal{A}(f)} \int_\Omega u(w)f_1 \, d\mu + (1-a) \sup_{w \in \mathcal{A}(f)} \int_\Omega u(w)f_2 \, d\mu \\
&\le a \sup_{w \in \mathcal{A}(f_1)} \int_\Omega u(w)f_1 \, d\mu + (1-a) \sup_{w \in \mathcal{A}(f_2)} \int_\Omega u(w)f_2 \, d\mu \\
&= aN_u(f_1) + (1-a)N_u(f_2),
\end{aligned}$$

as desired.  $\square$

REMARK 2.  The (*Arrow–Pratt*) *index of relative risk aversion* is defined in mathematical finance as

$$\mathrm{RA}_u(x) = -\frac{xu''(x)}{u'(x)} > 0 \qquad \text{for } x \in \mathbb{R}^+.$$

For a twice differentiable utility function $u$, it is clear that $\ln u^{-1}$ is convex if and only if $\mathrm{RA}_u(x) \ge 1$ for all $x \in (0, +\infty)$. Hence from Proposition 13, we can easily deduce that if $\mathrm{RA}_u \ge 1$, then the function $H_u$ is convex. In particular, this is true for the isoelastic utility $u_\gamma$ of order $\gamma < 0$ ($\mathrm{RA}_{u_\gamma} = 1 - \gamma$) and for the logarithmic utility $u_0$ ($\mathrm{RA}_{u_0} = 1$).

2.3. *Formula for $u$-entropy.*  Some results in this section can be deduced from [17]. They are presented here with complete proofs to make the present paper self-contained.

NOTATION 1. Let $\Lambda > 0$ and let $f \in D(\mu)$. We put $u^*(\Lambda/f)f = 0$ on $\{f = 0\}$, which is consistent with the limit $\lim_{x \to 0} u^*(\Lambda/x)x = 0$.

PROPOSITION 14. *For any $f \in D(\mu)$ and $\Lambda > 0$,*

$$N_u(f) \leq \int_\Omega u^*(\Lambda/f)f \, d\mu + \Lambda.$$

PROOF. Let $w \in \mathcal{A}_0(f)$. Then $u(w) \leq u^*(\Lambda/f) + (\Lambda/f)w$ on $\{f > 0\}$ because of (5). Multiplying by $f$ and integrating over $\{f > 0\}$ with respect to $\mu$, we get $\int_\Omega u(w)f \, d\mu \leq \int_\Omega u^*(\Lambda/f)f \, d\mu + \Lambda$. Then we take the supremum of the left-hand side over all $w \in \mathcal{A}_0(f)$ and apply Proposition 12 to obtain the assertion. □

PROPOSITION 15. *Let $f \in L^\infty(\mu)$. Then*

$$H_u(f) \leq \ln \|f\|_\infty.$$

PROOF. Put $K = \|f\|_\infty$. Let $\Lambda > 0$. Applying Proposition 14, we get

$$N_u(f) \leq \int_\Omega u^*(\Lambda/f)f \, d\mu + \Lambda \leq u^*(\Lambda/K) + \Lambda.$$

The last inequality holds because $u^*$ is a decreasing function. Taking the infimum of the right-hand side over all $\Lambda > 0$, we find that $N_u(f) \leq u(K)$ by (5). This implies that $H_u(f) \leq \ln K$. □

NOTATION 2. In this section we assume that an $f \in D(\mu)$ and a utility function $u : \mathbb{R}^+ \to \mathbb{R}$ are given. We define for any $x > 0$,

$$N_u(f; x) = \sup_{w \in \mathcal{A}(f;x)} \int_\Omega u(xw)f \, d\mu,$$

where

$$\mathcal{A}(f; x) = \{w \in D(\mu) : u(xw)^- \in L^1(f\mu)\}.$$

Clearly,

$$N_u(f) = N_u(f; 1) \quad \text{and} \quad \mathcal{A}(f) = \mathcal{A}(f; 1).$$

In the sequel we often write $N(x)$ for $N_u(f; x)$ and $\mathcal{A}(x)$ for $\mathcal{A}(f; x)$ if $u$ and $f$ are unambiguous.

PROPOSITION 16. *Let $f \in D(\mu)$ and let $H_u(f) < \infty$. Then:*

1. *The function $N : \mathbb{R}^+ \to \mathbb{R}$ is concave and increasing.*
2. *There exists a $y > 0$ such that $\sup_{x>0}(N(x) - yx) < \infty$.*

3. *If* $\mathrm{AE}(u) < 1$, *then* $\sup_{x>0}(N(x) - yx) < \infty$ *for every* $y > 0$.

PROOF. 1. First we show that $N$ is increasing. Let $0 < x_1 \leq x_2$. Since $u^-$ is decreasing, $\int_\Omega u(x_2 w)^- f\, d\mu \leq \int_\Omega u(x_1 w)^- f\, d\mu$. Therefore, $\mathcal{A}(x_1) \subset \mathcal{A}(x_2)$. Hence we deduce that

$$N(x_1) = \sup_{w \in \mathcal{A}(x_1)} \int_\Omega u(x_1 w) f\, d\mu$$

$$\leq \sup_{w \in \mathcal{A}(x_1)} \int_\Omega u(x_2 w) f\, d\mu$$

$$\leq \sup_{w \in \mathcal{A}(x_2)} \int_\Omega u(x_2 w) f\, d\mu = N(x_2).$$

Next we prove that $N$ is concave. Let $p \in (0,1)$, $x_1, x_2 > 0$. Put $x = px_1 + (1-p)x_2$. We first observe that $w_1 \in \mathcal{A}(x_1)$ and $w_2 \in \mathcal{A}(x_2)$ implies $w = x^{-1}(px_1 w_1 + (1-p)x_2 w_2) \in \mathcal{A}(x)$, since by the convexity of $u^-$,

$$\int_\Omega u(xw)^- f\, d\mu = \int_\Omega u(px_1 w_1 + (1-p)x_2 w_2)^- f\, d\mu$$

$$\leq \int_\Omega (pu(x_1 w_1)^- + (1-p)u(x_2 w_2)^-) f\, d\mu < \infty.$$

From this we conclude that

$$pN(x_1) + (1-p)N(x_2)$$

$$= p \sup_{w_1 \in \mathcal{A}(x_1)} \int_\Omega u(x_1 w_1) f\, d\mu + (1-p) \sup_{w_2 \in \mathcal{A}(x_2)} \int_\Omega u(x_2 w_2) f\, d\mu$$

$$= \sup_{\substack{w_1 \in \mathcal{A}(x_1) \\ w_2 \in \mathcal{A}(x_2)}} \int_\Omega (pu(x_1 w_1) + (1-p)u(x_2 w_2)) f\, d\mu$$

$$\leq \sup_{\substack{w_1 \in \mathcal{A}(x_1) \\ w_2 \in \mathcal{A}(x_2)}} \int_\Omega u(px_1 w_1 + (1-p)x_2 w_2) f\, d\mu$$

$$\leq \sup_{w \in \mathcal{A}(x)} \int_\Omega u(xw) f\, d\mu = N(x).$$

The last inequality holds because $x^{-1}(px_1 w_1 + (1-p)x_2 w_2) \in \mathcal{A}(x)$. Finally, we show that $-\infty < N(x) < \infty$ for any $x > 0$. Taking $w \equiv 1$, we can see that $N(x) \geq u(x) > -\infty$ for any $x > 0$. Since $N$ is increasing and $N(1) = N_u(f) < \infty$ by Proposition 9, we get $N(x) < \infty$ for any $x \leq 1$. On the other hand, for any $x > 1$ we take $p \in (0,1)$ such that $px + (1-p)x^{-1} = 1$ and deduce from the concavity of $N$ that $N(x) \leq p^{-1}(N(1) - (1-p)N(x^{-1})) < \infty$, which completes the proof.

2. Since $N$ is concave, it follows from Lemma 3 that there exist an $a > 0$ and a $b \in \mathbb{R}$ with the property $N(x) \le ax + b$ for all $x > 0$. Thus $\sup_{x>0}(N(x) - ax) \le b < \infty$.

3. First consider the case $u(\infty) > 0$. Let $\gamma = \mathrm{AE}(u) < 1$. By Proposition 7 there is an $x_0 > 0$ such that $0 < u(\lambda x) \le \lambda^\gamma u(x)$ for all $\lambda \ge 1$ and $x \ge x_0$. Then for any $\lambda \ge 1$ and $x > 0$,

$$
\begin{aligned}
N(\lambda x) &= \sup_{w \in \mathcal{A}(\lambda x)} \int_\Omega u(\lambda x w) f \, d\mu \\
&\le \sup_{w \in \mathcal{A}(\lambda x)} \int_\Omega u(\lambda (xw + x_0)) f \, d\mu \\
&\le \lambda^\gamma \sup_{w \in \mathcal{A}(\lambda x)} \int_\Omega u(xw + x_0) f \, d\mu \\
&\le \lambda^\gamma \sup_{v \in \mathcal{A}(x+x_0)} \int_\Omega u((x+x_0)v) f \, d\mu = \lambda^\gamma N(x + x_0).
\end{aligned}
$$

The last inequality holds because $0 < u(x_0) \le u(xw + x_0)$ for any $w \in \mathcal{A}(\lambda x)$ and consequently $(xw + x_0)/(x + x_0) \in \mathcal{A}(x + x_0)$.

From the above, for any $x > 0$,

$$
(6) \qquad N(x) \le N(x + x_0) \le \left( \frac{x + x_0}{x_0} \right)^\gamma N(2x_0) \le ax^\gamma + b,
$$

where

$$
\begin{aligned}
a &= x_0^{-\gamma} N(2x_0) \ge x_0^{-\gamma} \int_\Omega u(2x_0) f \, d\mu = x_0^{-\gamma} u(2x_0) > 0, \\
b &= N(2x_0).
\end{aligned}
$$

Now let $y > 0$. From (6) we deduce that there is a $c > 0$ such that $N(x) \le ax^\gamma + b \le \frac{1}{2} yx + c$ for all $x > 0$. This means that $\sup_{x>0}(N(x) - yx) \le \sup_{x>0}(-\frac{1}{2} xy + c) = c < \infty$, as required.

We now turn to the case $u(\infty) \le 0$. Let $\tilde{u} = u - u(\infty) + 1$. Then, according to Corollary 6, $\tilde{u}$ is a utility function and $\mathrm{AE}(\tilde{u}) < 1$. Moreover, $\tilde{u}(\infty) > 0$ and $N_u - u(\infty) + 1 = N_{\tilde{u}}$. Hence the proof of the assertion reduces to the first case. □

LEMMA 17.    *Let $f \in D(\mu)$ and let $\Lambda > 0$. Then*

$$
(7) \qquad \int_\Omega u^*(\Lambda / f) f \, d\mu \le \sup_{x>0}(N(x) - \Lambda x).
$$

Proof.    We put

$$u_n^*(y) = \sup_{n \geq x > 0} (u(x) - yx)$$

for any $n = 1, 2, \ldots$ and $y > 0$. Then $u_n^*(y) = u(I_n(y)) - yI_n(y)$, where $I_n(y) = \min\{I(y), n\}$. Moreover, $u(1) - y \leq u_n^*(y)$. Putting $y = \Lambda/f$, multiplying by $f$ and integrating over $\mu$, we obtain

$$
\begin{aligned}
(8) \quad u(1) - \Lambda &\leq \int_\Omega u_n^*(\Lambda/f) f \, d\mu \\
&= \int_\Omega (u(I_n(\Lambda/f)) - (\Lambda/f)I_n(\Lambda/f)) f \, d\mu \\
&= \int_\Omega u(x_n w_n) f \, d\mu - x_n \Lambda
\end{aligned}
$$

for any $\Lambda > 0$, where we put $I_n(\Lambda/f) = 0$ on $\{f = 0\}$ and where $0 < x_n = \int_\Omega I_n(\Lambda/f) \, d\mu \leq n < \infty$ and $w_n = x_n^{-1} I_n(\Lambda/f)$. Observe that $w_n \in \mathcal{A}(x_n)$, since $\int_\Omega u(x_n w_n)^- f \, d\mu \leq \Lambda - x_n \Lambda - u(1) + u(n)^+ < \infty$ by (8). As a result,

$$\int_\Omega u_n^*(\Lambda/f) f \, d\mu \leq \sup_{x > 0} (N(x) - \Lambda x).$$

Because $u(1)f - \Lambda \leq u_n^*(\Lambda/f)f \nearrow u^*(\Lambda/f)f$ pointwise on $\{f > 0\}$ as $n \to \infty$, by monotone convergence we obtain (7). □

Proposition 18.    *Assume that $f \in D(\mu)$ and $H_u(f) < \infty$. Then:*

1. *For every $\Lambda > 0$,*

$$-\infty < u^*(\Lambda) \leq \int_\Omega u^*(\Lambda/f) f \, d\mu.$$

2. *There exists a $\Lambda_0 > 0$ such that for each $\Lambda \geq \Lambda_0$,*

$$\int_\Omega u^*(\Lambda/f) f \, d\mu < \infty.$$

3. *If $\mathrm{AE}(u) < 1$, then for each $\Lambda > 0$*

$$\int_\Omega u^*(\Lambda/f) f \, d\mu < \infty.$$

Proof.    1. Let $\Lambda > 0$. Then, applying (3) for $y = \Lambda/f$, we get $u(x)f - \Lambda x \leq u^*(\Lambda/f)f$ on the set $\{f > 0\}$ for each $x > 0$. Hence, integrating with respect to $\mu$, taking the supremum over all $x > 0$ and applying (3) once again, we obtain $u^*(\Lambda) \leq \int_\Omega u^*(\Lambda/f) f \, d\mu$, as desired.

2. It is enough to prove that $\int_\Omega u^*(\Lambda_0/f) f \, d\mu < \infty$ for some $\Lambda_0 > 0$, because for $\Lambda \geq \Lambda_0$, we have $\int_\Omega u^*(\Lambda/f) f \, d\mu \leq \int_\Omega u^*(\Lambda_0/f) f \, d\mu < \infty$, since $u^*$ is nonincreasing. Now, the assertion follows from Proposition 16.2 and Lemma 17.

3. This follows from Proposition 16.3 and Lemma 17. □

PROPOSITION 19. *Let $f \in D(\mu)$. Then the following conditions are equivalent:*

1. *The inequality $H_u(f) < \infty$ holds.*
2. *There exists a $\Lambda_0 > 0$ such that*
$$\int_\Omega u^*(\Lambda_0/f) f \, d\mu < \infty.$$
3. *There exists a $\Lambda_0 > 0$ such that for each $\Lambda \geq \Lambda_0$*
$$\int_\Omega u^*(\Lambda/f) f \, d\mu < \infty.$$

PROOF. The implication $1 \Rightarrow 3$ follows from Proposition 18.2, $3 \Rightarrow 2$ is obvious and $2 \Rightarrow 1$ follows from Propositions 9 and 14. $\square$

NOTATION 3. Let $\Lambda > 0$ and let $f \in D(\mu)$. We put $I(\Lambda/f) = 0$ on $\{f = 0\}$.

THEOREM 20. *Let $f \in D(\mu)$. Assume that $\mathrm{AE}(u) < 1$ and $H_u(f) < \infty$. Then:*

1. *We have $\int_\Omega I(\Lambda/f) \, d\mu \in \mathbb{R}^+$ for all $\Lambda > 0$.*
2. *There exists a unique $\Lambda_f > 0$ such that*
$$\int_\Omega I(\Lambda_f/f) \, d\mu = 1. \tag{9}$$
3. *We have $I(\Lambda_f/f) \in \mathcal{A}_0(f)$.*
4. *The following formulae hold:*
$$N_u(f) = \int_\Omega u(I(\Lambda_f/f)) f \, d\mu = \int_\Omega u^*(\Lambda_f/f) f \, d\mu + \Lambda_f, \tag{10}$$
$$H_u(f) = \ln u^{-1} \left( \int_\Omega u(I(\Lambda_f/f)) f \, d\mu \right). \tag{11}$$

PROOF. 1. Since $u^*$ is convex, $I(\Lambda/f)(\Lambda - \widetilde{\Lambda}) = -(u^*)'(\Lambda/f)(\Lambda - \widetilde{\Lambda}) \leq (u^*(\widetilde{\Lambda}/f) - u^*(\Lambda/f)) f$ on $\{f > 0\}$ for any $\Lambda > \widetilde{\Lambda} > 0$. It follows that
$$(\Lambda - \widetilde{\Lambda}) \int_\Omega I(\Lambda/f) \, d\mu \leq \int_\Omega I(\Lambda/f)(\Lambda/f - \widetilde{\Lambda}/f) f \, d\mu$$
$$\leq \int_\Omega (u^*(\widetilde{\Lambda}/f) - u^*(\Lambda/f)) f \, d\mu < \infty,$$
because $-\infty < \int_\Omega u^*(\Lambda/f) f \, d\mu$ and $\int_\Omega u^*(\widetilde{\Lambda}/f) f \, d\mu < \infty$ by Proposition 18.1 and 18.3. As a result, $\int_\Omega I(\Lambda/f) \, d\mu < \infty$. Moreover, since $0 < I(\Lambda/f)$ on $\{f > 0\}$, it follows that $0 < \int_\Omega I(\Lambda/f) \, d\mu$, as desired.

2. Statement 1 above means that $\int_\Omega I(\Lambda/f)\,d\mu$ as a function of $\Lambda \in (0,\infty)$ has values in $(0,\infty)$. It is a strictly decreasing function because $I$ is. It is continuous with limit 0 as $\Lambda \nearrow \infty$ and $\infty$ as $\Lambda \searrow 0$ by monotone convergence, because $I$ has the same properties. As a result, there is a unique $\Lambda_f > 0$ such that $\int_\Omega I(\Lambda_f/f)\,d\mu = 1$.

3. Clearly, $I(\Lambda_f/f) \geq 0$, so it is a density by statement 2 above. Moreover, $I(\Lambda_f/f) = 0$ on $\{f = 0\}$ by the notational convention adopted (Notation 3).

4. Proposition 2 yields $u(x) - yx \leq u^*(y) = u(I(y)) - yI(y)$ for all $x, y > 0$. This implies that

$$u(w) - (\Lambda_f/f)w \leq u(I(\Lambda_f/f)) - (\Lambda_f/f)I(\Lambda_f/f) = u^*(\Lambda_f/f)$$

on $\{f > 0\}$ for any $w \in \mathcal{A}_0(f)$. Integrating with respect to $f\mu$, we get

$$\int_\Omega u(w)f\,d\mu \leq \int_\Omega u(I(\Lambda_f/f))f\,d\mu = \int_\Omega u^*(\Lambda_f/f)f\,d\mu + \Lambda_f$$

for any $w \in \mathcal{A}_0(f)$. Taking the supremum of the left-hand side over all such $w$'s and applying Proposition 12, we obtain (10) because $I(\Lambda_f/f) \in \mathcal{A}_0(f)$. Finally, (11) follows immediately from (10).  $\square$

REMARK 3.  To use (11) to compute the $u$-entropy $H_u(f)$, we need to know the constant $\Lambda_f$ defined implicitly by (9). Although the constant is determined uniquely by (9), it may not be possible to find a closed-form expression for it, except in some particular, though important, cases such as the logarithmic or isoelastic utility (see below). In general, the fact that $\Lambda_f$ is only defined implicitly is a limitation in using formula (11) for $H_u(f)$.

REMARK 4.  Adopting the methods of Kramkov and Schachermayer ([17], Section 5), it is possible to show that $\mathrm{AE}(u) < 1$ is a minimal assumption on the utility function $u$ for the validity of Theorem 20. If $\mathrm{AE}(u) = 1$, then the supremum

$$N_u(f) = \sup_{w \in \mathcal{A}(f)} \int_\Omega u(w)f\,d\mu$$

may fail to be attained at any $w \in \mathcal{A}(f)$, invalidating the assertions of Theorem 20. According to Proposition 18.3, the condition $\mathrm{AE}(u) < 1$ implies that

$$(12) \qquad \int_\Omega u^*(\Lambda/f)f\,d\mu < \infty \qquad \text{for all } \Lambda > 0.$$

By a similar argument as in [18], it can be demonstrated that (12) is in fact a necessary and sufficient condition for the assertions of Theorem 20 to hold.

2.4. *Examples.*

EXAMPLE 2 (Logarithmic utility). Let $u \colon \mathbb{R}^+ \to \mathbb{R}$ be given by $u(x) = \ln x$ for $x \in \mathbb{R}^+$. Then $H_u$ is equal to the *Boltzmann–Shannon conditional entropy* $H_1$ given by

$$H_1(f) = \int_\Omega f \ln f \, d\mu$$

for $f \in D(\mu)$.

EXAMPLE 3 (Isoelastic utility). Let $u \colon \mathbb{R}^+ \to \mathbb{R}$ be given by $u(x) = \frac{1}{\gamma} x^\gamma$ for $\gamma \in (-\infty, 0) \cup (0, 1)$ and $x \in \mathbb{R}^+$. Then $H_u$ is equal to the *Rényi entropy* $H_a$ of order $\alpha = (1 - \gamma)^{-1} \in (0, 1) \cup (1, \infty)$ given by

$$H_\alpha(f) = \frac{1}{\alpha - 1} \ln \int_\Omega f^\alpha \, d\mu = \frac{\alpha}{\alpha - 1} \ln \|f\|_\alpha$$

for any $f \in D(\mu) \cap L^\alpha(\mu)$.

**3. Markov operators.** Here we collect the definitions and properties that involve Markov operators. They can be found, for example, in [20].

DEFINITION 5. Let $P \colon D(\mu) \to D(\mu)$. We say that $P$ is a *Markov* (or *stochastic*) *operator on densities* if

(13) $$P(\lambda f_1 + (1 - \lambda) f_2) = \lambda P(f_1) + (1 - \lambda) P(f_2)$$

for all $f_1, f_2 \in D(\mu)$ and $\lambda \in [0, 1]$.

REMARK 5. A Markov operator $P$ on $D(\mu)$ can be uniquely extended to an operator $\bar{P} \colon L^1(\mu) \to L^1(\mu)$ such that:

1. $\bar{P}$ is linear;
2. $\bar{P} f \geq 0$ for every $0 \leq f \in L^1(\mu)$;
3. $\|\bar{P} f\|_1 = \|f\|_1$ for every $0 \leq f \in L^1(\mu)$.

The extended operator satisfies the condition:

4. $\|\bar{P} f\|_1 \leq \|f\|_1$ for every $f \in L^1(\mu)$ [i.e., $\bar{P}$ is a contraction on $L^1(\mu)$].

We call $\bar{P}$ a *Markov operator* in $L^1(\mu)$ or simply a Markov operator. For simplicity we use the same letter $P$ for a Markov operator on $D(\mu)$ and for its extension to $L^1(\mu)$.

PROPOSITION 21. *Let* $P \colon L^1(\mu) \to L^1(\mu)$ *be a Markov operator. Then there exists a unique operator* $P^* \colon L^\infty(\mu) \to L^\infty(\mu)$ *such that:*

1. $P^*$ is linear;
2. $P^*1 = 1$;
3. $P^*g \geq 0$ for every $0 \leq g \in L^\infty(\mu)$;
4. $\|P^*g\|_\infty \leq \|g\|_\infty$ for every $g \in L^\infty(\mu)$ [i.e., $P^*$ is a contraction on $L^\infty(\mu)$];
5. for all $f \in L^1(\mu)$ and $g \in L^\infty(\mu)$,

$$(14) \qquad \int_\Omega (Pf)g \, d\mu = \int_\Omega f(P^*g) \, d\mu.$$

We call $P^*$ the adjoint operator to $P$.

PROPOSITION 22 (Jensen inequality).   Let $P : L^1(\mu) \to L^1(\mu)$ be a Markov operator, let $u : \mathbb{R}^+ \to \mathbb{R}$ be a concave function and suppose that $g, u(g) \in L^\infty(\mu)$. Then

$$P^*(u(g)) \leq u(P^*g).$$

PROOF.   The idea of this proof is from [21]. Let $g \in L^\infty(\mu)$. According to Lemma 3, we can find an $a_n > 0$ and a $b_n \in \mathbb{R}$ for any $n \in \mathbb{N}$ such that $u(x) = \inf\{a_n x + b_n : n \in \mathbb{N}\}$. By Proposition 21,

$$\begin{aligned} P^*(u(g)) &= P^*(\inf\{a_n g + b_n : n \in \mathbb{N}\}) \\ &\leq \inf\{P^*(a_n g + b_n) : n \in \mathbb{N}\} \\ &= \inf\{a_n P^*g + b_n : n \in \mathbb{N}\} = u(P^*g), \end{aligned}$$

as desired.   □

DEFINITION 6.   Let $P : L^1(\mu) \to L^1(\mu)$ be a Markov operator and let $f \in D(\mu)$. We say that $P$ is doubly stochastic if $P1 = 1$, that is, 1 is a stationary density for $P$. Note that $P$ is doubly stochastic if and only if $\|P^*g\|_1 = \|g\|_1$ for all $g \in L^\infty(\mu)$, $g \geq 0$.

DEFINITION 7.   Let $P : L^1(\mu) \to L^1(\mu)$ be a doubly stochastic operator. We say that $P$ is exact (asymptotically stable) if $P^n f \xrightarrow{L^1} 1$ as $n \to \infty$ for every $f \in D(\mu)$.

## 4. Evolution of entropy.

### 4.1. H-theorem.

THEOREM 23.   Let $P : L^1(\mu) \to L^1(\mu)$ be a doubly stochastic operator and let $f \in D(\mu)$. Then

$$H_u(Pf) \leq H_u(f).$$

PROOF. If $H_u(f) = \infty$, then the assertion is obvious. Assume that $H_u(f) < \infty$. Take $w \in \mathcal{A}_{0b}(Pf)$. Define $w_n = w \vee (1/n)$ and $x_n = \int_\Omega w_n \, d\mu$ for $n = 1$, $2, \ldots$. Then $w_n$ and $u(w_n)$ are bounded, $1 \le x_n \le 1 + 1/n$ and $0 < x_n^{-1} w_n \in D(\mu)$. Furthermore, $w_n \searrow w$, $u(w_n) \searrow u(w)$ and $x_n \searrow 1$ as $n \to \infty$. Moreover, $\int_\Omega u(w_n)^- Pf \, d\mu \le \int_\Omega u(w)^- Pf \, d\mu < \infty$ for $n \in \mathbb{N}$. Hence, applying the monotone convergence theorem, formula (14) and Proposition 22, we get

$$
\begin{aligned}
\int_\Omega u(w) Pf \, d\mu &= \lim_{n \to \infty} \int_\Omega u(w_n) Pf \, d\mu \\
&= \lim_{n \to \infty} \int_\Omega P^*(u(w_n)) f \, d\mu \\
&\le \lim_{n \to \infty} \int_\Omega u(P^* w_n) f \, d\mu \\
&= \lim_{n \to \infty} \int_\Omega u(x_n P^*(x_n^{-1} w_n)) f \, d\mu.
\end{aligned}
$$
(15)

For any $n \in \mathbb{N}$, we have

$$
\int_\Omega P^*(x_n^{-1} w_n) \, d\mu = x_n^{-1} \int w_n \, P1 \, d\mu = x_n^{-1} \int_\Omega w_n \, d\mu = 1.
$$
(16)

Using Proposition 22 and (14) once again (for the concave function $-u^-$), we get

$$
\begin{aligned}
\int_\Omega u(x_n P^*(x_n^{-1} w_n))^- f \, d\mu &= \int_\Omega u(P^*(w_n))^- f \, d\mu \\
&\le \int_\Omega P^*(u(w_n)^-) f \, d\mu \\
&= \int_\Omega u(w_n)^- Pf \, d\mu \\
&\le \int_\Omega u(w)^- Pf \, d\mu < \infty.
\end{aligned}
$$
(17)

From (16) and (17) we deduce that $P^*(x_n^{-1} w_n) \in \mathcal{A}(x_n)$ (see Notation 2). From this and (15), it follows that

$$
\begin{aligned}
\int_\Omega u(w) Pf \, d\mu &\le \lim_{n \to \infty} \int_\Omega u(x_n P^*(x_n^{-1} w_n)) f \, d\mu \\
&\le \lim_{n \to \infty} N(x_n).
\end{aligned}
$$

According to Proposition 16.1, $N$ is concave and hence continuous, so it follows that $\lim_{n \to \infty} N(x_n) = N(1)$. Thus

$$
\int_\Omega u(w) Pf \, d\mu \le N(1) = N_u(f).
$$

Taking the supremum over all $w \in \mathcal{A}_{0b}(Pf)$, we get $N_u(Pf) \leq N_u(f)$. The assertion of the theorem follows because $\ln \circ u^{-1}$ is an increasing function. $\square$

REMARK 6. Let $P \colon L^1(\mu) \to L^1(\mu)$ be a doubly stochastic operator and let $f \in D(\mu)$ be such that $H_u(f) < \infty$. As a consequence of Theorem 23, the sequence $H_u(P^n f)$, $n = 1, 2, \ldots$, is nonincreasing. If $P$ is invertible, then $H_u(P^n f)$, $n = 1, 2, \ldots$, is a constant sequence.

### 4.2. *Inequalities.*

LEMMA 24. *Suppose that $f \in D(\mu)$ and $w \in \mathcal{A}_{0b}(f)$. Then*

$$\int_\Omega u(w) f \, d\mu \leq u'(1) \|w\|_\infty \|f - 1\|_1 + u(1).$$

PROOF. Let $f \in D(\mu)$ and let $w \in \mathcal{A}_{0b}(f)$. By the concavity of $u$ we have $u(x) \leq (x-1)u'(1) + u(1)$ for any $x \geq 0$. Hence

$$\int_\Omega u(w) f \, d\mu \leq u'(1) \int_\Omega (w-1) f \, d\mu + u(1)$$

$$= u'(1) \int_\Omega (f-1) w \, d\mu + u(1) \leq u'(1) \|w\|_\infty \|f - 1\|_1 + u(1),$$

as required. $\square$

PROPOSITION 25. *Let $u \colon \mathbb{R}^+ \to \mathbb{R}$ be a utility function such that $\mathrm{AE}(u) < 1$ and let $f \in D(\mu) \cap L^\infty(\mu)$. Put $K = \|f\|_\infty \geq 1$. Then for each $0 < C < 1$,*

$$N_u(f) \leq u'(1) I\left(u'\left(\frac{K-C}{1-C}\right)\frac{C}{K}\right) \|f - 1\|_1 + u(1).$$

PROOF. Since

$$1 = \int_{\{f \geq C\}} f \, d\mu + \int_{\{f < C\}} f \, d\mu \leq K\mu\{f \geq C\} + C(1 - \mu\{f \geq C\}),$$

it follows that $\mu\{f \geq C\} \geq \frac{1-C}{K-C}$. Let us take $\Lambda_f$ as in Theorem 20. Then

$$1 = \int_\Omega I\left(\frac{\Lambda_f}{f}\right) d\mu \geq I\left(\frac{\Lambda_f}{C}\right)\mu\{f \geq C\} \geq I\left(\frac{\Lambda_f}{C}\right)\frac{1-C}{K-C}.$$

Consequently, $\Lambda_f \geq u'(\frac{K-C}{1-C})C$ and so

$$I\left(\frac{\Lambda_f}{f}\right) \leq I\left(u'\left(\frac{K-C}{1-C}\right)\frac{C}{f}\right) \leq I\left(u'\left(\frac{K-C}{1-C}\right)\frac{C}{K}\right).$$

Now, using Theorem 20 and Lemma 24, we have $I(\Lambda_f/f) \in \mathcal{A}_{0b}(f)$ and

$$N_u(f) = \int_\Omega u\left(I\left(\frac{\Lambda_f}{f}\right)\right) f \, d\mu \leq u'(1) I\left(u'\left(\frac{K-C}{1-C}\right)\frac{C}{K}\right) \|f - 1\|_1 + u(1),$$

as required. $\square$

4.3. *Main theorems.* Throughout this section we assume that $u \colon \mathbb{R}^+ \to \mathbb{R}$ is a utility function and that $P \colon L^1(\mu) \to L^1(\mu)$ is a doubly stochastic operator.

THEOREM 26. *Let us assume that* $\mathrm{AE}(u) < 1$. *Let* $f \in D(\mu)$ *be such that* $H_u(f) < \infty$. *Then*

$$P^n f \xrightarrow{L^1} 1 \quad \textit{as } n \to \infty \quad \implies \quad H_u(P^n f) \searrow 0 \quad \textit{as } n \to \infty.$$

PROOF. *Step* 1. First observe that the condition $H_u(P^n f) \to 0$ as $n \to \infty$ is equivalent to $N_u(P^n f) \to u(1)$ as $n \to \infty$. Moreover, according to Theorem 23, the sequence $N_u(P^n f)$ is nonincreasing.

*Step* 2. Let $f \in D(\mu) \cap L^\infty(\mu)$. Then the assertion follows immediately from Proposition 25 and the fact that $u(1) \le N_u(f)$.

*Step* 3. Now we assume that $f \in D(\mu) \setminus L^\infty(\mu)$. Then we can define two densities $f_c = (f/a_c)\mathbb{1}_{\{f < c\}}$ and $f^c = (f/a^c)\mathbb{1}_{\{f \ge c\}}$ for any $c > (\mu\{f > 0\})^{-1}$, where $a_c = \int_{\{f < c\}} f \, d\mu > 0$ and $a^c = \int_{\{f \ge c\}} f \, d\mu > 0$. Clearly, $a_c + a^c = 1$, $f = a_c f_c + a^c f^c$ and $P^n f = a_c P^n f_c + a^c P^n f^c$ for every $n = 1, 2, \ldots$ . Moreover, $f_c \in D(\mu) \cap L^\infty(\mu)$. According to Proposition 18.3, $\int_\Omega u^*(y/f) f \, d\mu < \infty$ for each $y > 0$. Let $\Lambda > 0$. Then $\int_\Omega u^*(\Lambda/f^c) f^c \, d\mu = (1/a^c) \int_{\{f \ge c\}} u^*(\Lambda a^c/f) f \, d\mu < \infty$. Applying Proposition 19, we deduce that $H_u(f^c) < \infty$. According to Proposition 13,

$$(18) \qquad N_u(P^n f) \le a_c N_u(P^n f_c) + a^c N_u(P^n f^c)$$

for each $n \in \mathbb{N}$. To estimate the second term, observe that by Theorems 23 and 20 there is a $w^c \in \mathcal{A}_0(f^c)$ such that

$$a^c N_u(P^n f^c) \le a^c N_u(f^c) = a^c \int_\Omega u(w^c) f^c \, d\mu = \int_{\{f \ge c\}} u(w^c) f \, d\mu.$$

Let $\varepsilon > 0$. Since $u(w^c) \le u^*(\frac{\varepsilon}{3f}) + \frac{\varepsilon}{3f} w^c$ on $\{f \ge c\}$ by (3),

$$(19) \qquad a^c N_u(P^n f^c) \le \int_{\{f \ge c\}} u^*\left(\frac{\varepsilon}{3f}\right) f \, d\mu + \frac{\varepsilon}{3}.$$

As $\int_\Omega u^*(\frac{\varepsilon}{3f}) f \, d\mu < \infty$ by Proposition 18.3, there is a $c > (\mu\{f > 0\})^{-1}$ such that

$$(20) \qquad \int_{\{f \ge c\}} u^*\left(\frac{\varepsilon}{3f}\right) f \, d\mu \le \frac{\varepsilon}{3}.$$

It follows from Step 2 that the first term on the right-hand side of (18) tends to $a_c u(1)$. We can therefore take an $N \in \mathbb{N}$ such that for each $n > N$,

$$(21) \qquad a_c N_u(P^n f_c) \le u(1) + \varepsilon/3.$$

Inequalities (18)–(21) and Proposition 8 lead to

$$u(1) \leq N_u(P^n f) \leq u(1) + \varepsilon$$

for each $n > N$, which completes the proof. $\square$

THEOREM 27.  *Let $f \in D(\mu)$. Then*

$$H_u(P^n f) \to 0 \quad \text{as } n \to \infty \quad \Longrightarrow \quad P^n f \xrightarrow{L^1} 1 \quad \text{as } n \to \infty.$$

PROOF.  By Proposition 11 the theorem is a straightforward consequence of the lemma below. $\square$

LEMMA 28.  *Let $u(1) = 0$ and $f_n \in D(\mu)$ for $n = 1, 2, \ldots$. Then*

$$N_u(f_n) \to 0 \quad \text{as } n \to \infty \quad \Longrightarrow \quad f_n \xrightarrow{L^1} 1 \quad \text{as } n \to \infty.$$

PROOF.  Clearly, we can assume that none of the $f_n$'s is identically ($\mu$-a.e.) equal to 1. We define $A_n = \{f_n \geq 1\}$ and $\Sigma_n = \{\varnothing, A_n, A_n^\urcorner, \Omega\}$. Put $0 < q_n = \mu(A_n) < 1$ and $p_n = \int_{A_n} f_n \, d\mu$ for $n = 1, 2, \ldots$. Then

$$E_\mu(f_n | \Sigma_n) = \frac{p_n}{q_n} \mathbb{1}_{A_n} + \frac{1 - p_n}{1 - q_n} \mathbb{1}_{A_n^\urcorner}.$$

Hence

$$(22) \qquad N_u(E_\mu(f_n | \Sigma_n)) = N(p_n, q_n),$$

where

$$N(p, q) = \sup\{u(w_1)p + u(w_2)(1 - p) : w_1 q + w_2(1 - q) = 1, w_1, w_2 \geq 0\}$$

for $0 \leq q \leq p \leq 1$. Because $L^1(\mu) \ni g \mapsto E_\mu(g | \Sigma_n) \in L^1(\mu)$ is a doubly stochastic operator, we can apply Theorem 23 to get

$$(23) \qquad N_u(f_n) \geq N_u(E_\mu(f_n | \Sigma_n)).$$

From (22) and (23) it follows that $N(p_n, q_n) \to 0$ as $n \to \infty$.

Moreover, it is easy to check that

$$\|f_n - 1\|_1 = \|E_\mu(f_n | \Sigma_n) - 1\|_1 = 2|p_n - q_n|.$$

Suppose, contrary to our claim, that $\|f_n - 1\|_1 \nrightarrow 0$ as $n \to \infty$. This would mean that $|p_n - q_n| \nrightarrow 0$ as $n \to \infty$. Therefore, by passing to a subsequence if necessary, we can assume that $(p_n, q_n) \to (p, q)$ such that $p \neq q$. Then $N(p, q) > 0$ and hence there exist $w_1, w_2 \geq 0$ such that $w_1 q + w_2(1 - q) = 1$ and $\delta = u(w_1)p + u(w_2)(1 - p) > 0$. Let us consider two cases.

*Case 1.* $u(w_1), u(w_2) \neq -\infty$. Then $0 < z_n = w_1 q_n + w_2(1 - q_n) \to 1$ as $n \to \infty$. Put $w_1^n = w_1 / z_n$ and $w_2^n = w_2 / z_n$ for $n = 1, 2, \ldots$. We have $w_1^n q_n + w_2^n(1 -$ $q_n) = 1$. Hence $N(p_n, q_n) \geq u(w_1^n)p_n + u(w_2^n)(1 - p_n) \to \delta > 0$ as $n \to \infty$, a contradiction.

*Case* 2. $u(w_1) = -\infty$ or $u(w_2) = -\infty$. We can assume that $u(w_1) = -\infty$, since the other case is similar. Then $w_1 = 0$, $w_2 = \frac{1}{1-q} > 1$, $p = 0$ and $q < 1$. Put $w_1^n = u^{-1}(-1/\sqrt{p_n})$ and $w_2^n = (1 - w_1^n q_n)/(1 - q_n)$ for $n = 1, 2, \ldots$. Then $w_1^n q_n + w_2^n(1 - q_n) = 1$, $w_1^n \to 0$ and $w_2^n \to \frac{1}{1-q}$ as $n \to \infty$. Hence $N(p_n, q_n) \geq u(w_1^n)p_n + u(w_2^n)(1 - p_n) \to u(\frac{1}{1-q}) = \delta > 0$ as $n \to \infty$, a contradiction. $\square$

THEOREM 29. *Suppose that* $\mathrm{AE}(u) < 1$. *Then the following conditions are equivalent:*

1. *$P$ is exact;*
2. *$H_u(P^n f) \to 0$ as $n \to \infty$ for all $f \in L^\infty(\mu) \cap D(\mu)$;*
3. *$H_u(P^n f) \to 0$ as $n \to \infty$ for all $f \in D(\mu)$ such that $H_u(f) < \infty$.*

PROOF. $1 \Rightarrow 3$. The assertion follows from Theorem 26.

$3 \Rightarrow 2$. Obvious.

$2 \Rightarrow 1$. By Theorem 27, $P^n f \xrightarrow{L^1} 1$ as $n \to \infty$ for all $f \in L^\infty(\mu) \cap D(\mu)$. Now let $f \in D(\mu)$ be an arbitrary density. Take any $\varepsilon > 0$. There exists a $c > 0$ such that $\int_{\{f \geq c\}} f \, d\mu \leq \varepsilon/4$. Let $f_c = (\int_{\{f < c\}} f \, d\mu)^{-1} \mathbb{1}_{\{f < c\}} f$. Then $f_c \in L^\infty(\mu) \cap D(\mu)$, and so we can find $N \in \mathbb{N}$ such that $\|P^n f_c - 1\|_1 \leq \varepsilon/2$ for every $n \geq N$. Finally, we get

$$\|P^n f - 1\|_1 \leq \|P^n f - P^n f_c\|_1 + \|P^n f_c - 1\|_1$$

$$\leq \|f - f_c\|_1 + \varepsilon/2 = 2 \int_{f \geq c} f \, d\mu + \varepsilon/2 \leq \varepsilon,$$

which completes the proof. $\square$

COROLLARY 30. *Let* $\alpha \in (0, 1) \cup (1, +\infty)$. *The following conditions are equivalent:*

1. *$P$ is exact;*
2. *$\|P^n f\|_\alpha \to 1$ as $n \to \infty$ for all $f \in L^\infty(\mu) \cap D(\mu)$;*
3. *$\|P^n f\|_\alpha \to 1$ as $n \to \infty$ for all $f \in L^\alpha(\mu) \cap D(\mu)$.*

PROOF. This follows because $H_\alpha(f) = \frac{\alpha}{\alpha-1} \ln \|f\|_\alpha$; see Example 3. $\square$

REMARK 7. The implication $3 \Rightarrow 1$ in Corollary 30 follows from a result proved by Loskot and Rudnicki [21], already mentioned in the Introduction. The reverse implication $1 \Rightarrow 3$ for $0 < \alpha < 1$ is an easy exercise involving integral inequalities.

### 4.4. *Continuous time.*

DEFINITION 8. We call a family of operators $\{P_t\}_{t\geq 0}$ a *doubly stochastic semigroup* if:

1. $P_t : L^1(\mu) \to L^1(\mu)$ is a doubly stochastic operator for each $t \geq 0$;
2. $P_t \circ P_s = P_{t+s}$ for any $t, s \geq 0$;
3. $P_0 = \mathrm{Id}_\Omega$.

We say that $\{P_t\}_{t\geq 0}$ is *asymptotically stable* if $P^t f \xrightarrow{L^1} 1$ as $t \to \infty$ for every $f \in D(\mu)$.

THEOREM 31. *Suppose that* $\mathrm{AE}(u) < 1$. *Then the following conditions are equivalent:*

1. $\{P_t\}_{t\geq 0}$ *is asymptotically stable;*
2. $H_u(P^t f) \to 0$ *as* $t \to \infty$ *for each* $f \in D(\mu)$ *such that* $H_u(f) < \infty$.

PROOF. $1 \Rightarrow 2$. Let $f \in D(\mu)$. According to Theorem 23, the function $[0, \infty) \ni t \to H_u(P^t f) \in \mathbb{R}$ is decreasing. Moreover, it follows from Theorem 29 that $H_u(P^n f) \to 0$ as $n \to \infty$. These two statements imply the assertion.

$2 \Rightarrow 1$. The assertion follows from Proposition 11 and Lemma 28. $\quad\square$

Institute of Mathematics
Jagiellonian University
Reymonta 4
30-059 Kraków
Poland
E-mail: wojciech.slomczynski@im.uj.edu.pl

Department of Mathematics
University of Hull
Cottingham Road
Kingston upon Hull HU6 7RX
United Kingdom
E-mail: t.j.zastawniak@hull.ac.uk